\documentclass[a4paper,twoside,reqno,12pt,tbtags]{article}

\usepackage{amssymb,amsmath,amsthm,mathrsfs}
\usepackage{graphicx}
\usepackage{times}

\numberwithin{equation}{section}

\newtheorem{prop}{Proposition}[section]
\newtheorem{thm}[prop]{Theorem}
\newtheorem{lem}[prop]{Lemma}
\newtheorem{cor}[prop]{Corollary}
\theoremstyle{definition}

\DeclareMathOperator{\conv}{conv}
\DeclareMathOperator{\interior}{int}

\newcommand{\EE}{\mathbb{E}}
\newcommand{\PP}{\mathbb{P}}
\newcommand{\RR}{\mathbb{R}}

\newcommand{\Rd}{\RR^d}
\newcommand{\compo}{\circ \hspace*{0.12em}}

\begin{document}

\thispagestyle{empty}

\begin{center}
	{\large University of Bern}\\
	{\large Institute of Mathematical Statistics and Actuarial Science}
	
	\textbf{\large Technical Report 74}
\end{center}
\strut

\vfill

\begin{center}
	\textbf{\Large Multivariate Log-Concave Distributions\\
	as a Nearly Parametric Model$^*$}

	{\large Dominic Schuhmacher, Andr\'{e} H\"usler and Lutz D\"umbgen
	
	July 2009 \ (minor revisions in February 2010 and April 2011)}
\end{center}

\vfill

\begin{abstract}
In this paper we show that the family $\mathcal{P}_d^{(\mathrm{lc})}$ of probability distributions on $\Rd$ with log-concave densities satisfies a strong continuity condition. In particular, it turns out that weak convergence within this family entails (i) convergence in total variation distance, (ii) convergence of arbitrary moments, and (iii) pointwise convergence of Laplace transforms. In this and several other respects the nonparametric model $\mathcal{P}_d^{(\mathrm{lc})}$ behaves like a parametric model such as, for instance, the family of all $d$-variate Gaussian distributions. As a consequence of the continuity result, we prove the existence of nontrivial confidence sets for the moments of an unknown distribution in $\mathcal{P}_d^{(\mathrm{lc})}$. Our results are based on various new inequalities for log-concave distributions which are of independent interest.
\end{abstract}

\vfill
\vfill



\noindent${}^*$ Work supported by Swiss National Science Foundation

\newpage

\section{Introduction}

It is well-known that certain statistical functionals such as moments fail to be weakly continuous on the set of, say, all probability measures on the real line for which these functionals are well-defined. This is the intrinsic reason why it is impossible to construct nontrivial two-sided confidence intervals for such functionals. For the mean and other moments, this fact was pointed out by Bahadur and Savage (1956). Donoho (1988) extended these considerations by noting that some functionals of interest are at least weakly semi-continuous, so that one-sided confidence bounds are possible. 

When looking at the proofs of the results just mentioned, one realizes that they often involve rather strange, e.g.\ multimodal or heavy-tailed, distributions. Natural questions are whether statistical functionals such as moments become weakly continuous and whether honest confidence intervals exist for these functionals if attention is restricted to a suitable nonparametric class of distributions. For instance, one possibility would be to focus on distributions on a given bounded region. But this may be too restrictive or lead to rather conservative procedures.

Alternatively we propose a qualitative constraint. When asking a statistician to draw a typical probability density, she or he will often sketch a bell-shaped, maybe skewed density. This suggests unimodality as a constraint, but this would not rule out heavy tails. In the present paper we favor the stronger though natural constraint of log-concavity, also called strong unimodality. One should note here that additional assumptions such as given bounded support or log-concavity can never be strictly verified based on empirical data alone; see Donoho~(1988, Section~2).

Before proceeding with log-concavity, let us consider briefly the parametric model $\mathcal{N}_d$ of all nondegenerate Gaussian distributions on $\Rd$. Suppose that a sequence of distributions $P_n = N_d(\mu_n, \Sigma_n) \in \mathcal{N}_d$ converges weakly to $P = N_d(\mu,\Sigma) \in \mathcal{N}_d$. This is easily shown to be equivalent to $\mu_n \to \mu$ and $\Sigma_n \to \Sigma$ as $n \to \infty$. But this implies convergence in total variation distance, i.e.
$$
	\lim_{n \to \infty} \int_{\Rd} |f_n(x) - f(x)| \, dx \ = \ 0 ,
$$
where $f_n$ and $f$ denote the Lebesgue densities of $P_n$ and $P$, respectively. Furthermore, weak convergence of $(P_n)_n$ to $P$ in $\mathcal{N}_d$ implies convergence of all moments and pointwise convergence of the Laplace-transforms. That means, for all $d$-variate polynomials $\Pi : \Rd \to \RR$,
$$
	\lim_{n \to \infty} \int \Pi(x) f_n(x) \, dx
	\ = \ \int \Pi(x) f(x) \, dx ,
$$
and for arbitrary $\theta \in \Rd$,
$$
	\lim_{n \to \infty} \int \exp(\theta^\top x) f_n(x) \, dx
	\ = \ \int \exp(\theta^\top x) f(x) \, dx .
$$

In the present paper we show that the nonparametric model $\mathcal{P}_d^{(\mathrm{lc})}$ of all log-concave probability distributions $P$ on $\Rd$ has the same properties. Log-concavity of $P$ means that it admits a Lebesgue density $f$ of the form
$$
	f(x) \ = \ \exp(\varphi(x))
$$
for some concave function $\varphi : \Rd \to [-\infty,\infty)$. Obviously the model $\mathcal{P}_d^{(\mathrm{lc})}$ contains the parametric family $\mathcal{N}_d$. All of its members are unimodal in that the level sets $\{x \in \Rd : f(x) \ge c\}$, $c > 0$, are bounded and convex. It is further known that product measures, marginals, convolutions, and weak limits (if a limiting density exists) of log-concave distributions are log-concave; see Dharmadhikari and Joag-dev (1988), Chapter~2. These closedness properties are again shared by the class of Gaussian distributions. The results in the present paper make a substantial contribution to the list of such shared properties and thus promote the view of the model $\mathcal{P}_d^{(\mathrm{lc})}$ as a viable nonparametric substitute for the Gaussian model $\mathcal{N}_d$.

The univariate class $\mathcal{P}_1^{(\mathrm{lc})}$ has been studied extensively; see Bagnoli and Bergstrom (2005), D\"umbgen and Rufibach (2009) and the references therein. Many standard models of univariate distributions belong to this nonparametric family, e.g.\ all gamma distributions with shape parameter $\ge 1$, and all beta distributions with both parameters $\ge 1$. Bagnoli and Bergstrom (2005) establish various properties of the corresponding distribution and hazard functions. Nonparametric maximum likelihood estimation of a distribution in $\mathcal{P}_1^{(\mathrm{lc})}$ has been studied by Pal et al.\ (2006) and D\"umbgen and Rufibach (2009). In particular, the latter two papers provide consistency results for these estimators. The findings of the present paper allow to strengthen these results considerably by showing that consistency in any reasonable sense implies consistency of all moments and, much more generally, consistency of the densities in exponentially weighted total variation distance. Algorithms for the one-dimensional maximum-likelihood estimator are described by D\"{u}mbgen et al.\ (2007) and D\"{u}mbgen and Rufibach~(2011).

The multivariate class $\mathcal{P}_d^{(\mathrm{lc})}$ is in various respects more difficult to treat. It has been considered in Dharmadhikari and Joag-dev (1988) and An~(1998). Comprehensive treatments of the state of the art in multivariate log-concave density modeling and estimation are Cule et al.\ (2010) and the survey paper by Walther (2009). An explicit algorithm for the nonparametric maximum likelihood estimator is provided by Cule et al.\ (2009). Consistency of this estimator has been verified by Cule and Samworth (2010) and Schuhmacher and D\"umbgen (2010). Again the results of the present paper allow to transfer consistency properties into much stronger modes of consistency.

The remainder of this paper is organized as follows. In Section~\ref{sec: Main results} we present our main result and some consequences, including an existence proof of non-trivial confidence sets for moments of log-concave distributions. Section~\ref{sec: Various inequalities} collects some basic inequalities for log-concave distributions which are essential for the main results and of independent interest. Most proofs are deferred to Section~\ref{sec: Proofs}.

\section{The main results}
\label{sec: Main results}

Let us first introduce some notation. Throughout this paper, $\|\cdot\|$ stands for Euclidean norm. The closed Euclidean ball with center $x \in \Rd$ and radius $\epsilon \ge 0$ is denoted by $B(x,\epsilon)$. With $\interior(S)$ and $\partial S$ we denote the interior and boundary, respectively, of a set $S \subset \Rd$.

\begin{thm}
\label{thm:dabigone}
Let $P$, $P_1$, $P_2$, $P_3$ \ldots be probability measures in $\mathcal{P}_d^{(\mathrm{lc})}$ with densities $f$, $f_1$, $f_2$, $f_3$, \ldots, respectively, such that $P_n \to P$ weakly as $n \to \infty$. Then the following two conclusions hold true:

\noindent
\textbf{(i)} \ The sequence $(f_n)$ converges uniformly to $f$ on any closed set of continuity points of $f$.

\noindent
\textbf{(ii)} \ Let $A: \Rd \to \RR$ be a sublinear function, i.e.\ $A(x+y) \le A(x) + A(y)$ and $A(rx) = rA(x)$ for all $x,y \in \Rd$ and $r \geq 0$. If
\begin{equation}
\label{eq:dacondition}
	f(x) \exp (A(x)) \ \to \ 0 \quad \text{as} \ \|x\| \to \infty,
\end{equation}
then $\int_{\Rd} \exp(A(x)) f(x) \, dx < \infty$ and
\begin{equation}
\label{eq:daintconvergence}
	\lim_{n \to \infty} \int_{\Rd} \exp(A(x)) \bigl| f_n(x) - f(x) \bigr| \, dx
	\ = \ 0 .
\end{equation}
\end{thm}

It is well-known from convex analysis that $\varphi = \log f$ is continuous on $\interior(\{\varphi > - \infty\}) = \interior(\{f > 0\})$. Hence the discontinuity points of $f$, if any, are contained in $\partial \{f > 0\}$. But $\{f > 0\}$ is a convex set, so its boundary has Lebesgue measure zero (cf.\ Lang 1986\nocite{lang86}). Therefore Part~(i) of Theorem~\ref{thm:dabigone} implies that $(f_n)_n$ converges to $f$ pointwise almost everywhere.

Note also that $f(x) \le C_1 \exp(- C_2 \|x\|)$ for suitable constants $C_1 = C_1(f) > 0$ and $C_2 = C_2(f) > 0$; see Corollary~\ref{cor:lutz} in Section~\ref{sec: Various inequalities}. Hence one may take $A(x) = c \|x\|$ for any $c \in [0, C_2)$ in order to satisfy (\ref{eq:dacondition}). Theorem~\ref{thm:dabigone} is a multivariate version of H\"{u}sler (2008, Theorem~2.1). It is also more general than findings of Cule and Samworth (2010) who treated the special case of $A(x) = \epsilon \|x\|$ for some small $\epsilon > 0$ with different techniques.

Before presenting the conclusions about moments and moment generating functions announced in the introduction, let us provide some information about the moment generating functions of distributions in $\mathcal{P}_d^{(\mathrm{lc})}$:

\begin{prop}
\label{prop:godot}
For a distribution $P \in \mathcal{P}_d^{(\mathrm{lc})}$ let $\Theta(P)$ be the set of all $\theta \in \Rd$ such that $\int \exp(\theta^\top x) \, P(dx) < \infty$. This set $\Theta(P)$ is convex, open and contains $0$. Let $\theta \in \Rd$ and $\epsilon > 0$ such that $B(\theta,\epsilon) \subset \Theta(P)$. Then
$$
	A(x) \ := \ \theta^\top x + \epsilon \|x\|
$$
defines a sublinear function $A$ on $\Rd$ such that the density $f$ of $P$ satisfies
$$
	\lim_{\|x\| \to \infty} \exp(A(x)) f(x) \ = \ 0 .
$$
\end{prop}

Note that for any $d$-variate polynomial $\Pi$ and arbitrary $\epsilon > 0$ there exists an $R = R(\Pi,\epsilon) > 0$ such that $|\Pi(x)| \le \exp(\epsilon \|x\|)$ for $\|x\|>R$. Hence part~(ii) of Theorem~\ref{thm:dabigone} and Proposition~\ref{prop:godot} entail the first part of the following theorem:

\begin{thm}
\label{thm:aebe}
Under the conditions of Theorem~\ref{thm:dabigone}, for any $\theta \in \Theta(P)$ and arbitrary $d$-variate polynomials $\Pi : \Rd \to \RR$, the integral $\int_{\Rd} \exp(\theta^\top x) |\Pi(x)| f(x) \, dx$ is finite and
$$
	\lim_{n \to \infty} \int_{\Rd} \exp(\theta^\top x) |\Pi(x)| \bigl| f_n(x) - f(x) \bigr| \, dx
	\ = \ 0 .
$$
Moreover, for any $\theta \in \Rd \setminus \Theta(P)$,
$$
	\lim_{n \to \infty} \int_{\Rd} \exp(\theta^\top x) f_n(x) \, dx
	\ = \ \infty .
$$
\end{thm}

\paragraph{Existence of nontrivial confidence sets for moments.}
With the previous results we can prove the existence of confidence sets for arbitrary moments, modifying Donoho's (1988) \nocite{Donoho_1988} recipe. Let $\mathcal{H} = \mathcal{H}_d$ denote the set of all closed halfspaces in $\Rd$. For two probability measures $P$ and $Q$ on $\Rd$ let
$$
	\|P - Q\|_{\mathcal{H}} \ := \ \sup_{H \in \mathcal{H}} \bigl| P(H) - Q(H) \bigr| .
$$
It is well-known from empirical process theory (e.g.\ van der Vaart and Wellner 1996\nocite{VanDerVaart_Wellner_1996}, Section 2.19) that for any $\alpha \in (0,1)$ there exists a universal constant $c_{\alpha,d}$ such that
$$
	\PP \Bigl( \bigl\| \hat{P}_n - P \bigr\|_{\mathcal{H}} \ge n^{-1/2} c_{\alpha,d} \Bigr)
	\ \le \ \alpha
$$
for arbitrary distributions $P$ on $\Rd$ and the empirical distribution $\hat{P}_n$ of independent random vectors $X_1, X_2, \ldots, X_n \sim P$. In particular, Massart's (1990) \nocite{Massart_1990} inequality yields the constant $c_{\alpha,1} = \bigl( \log(2/\alpha)/2 \bigr)^{1/2}$.

Under the assumption that $P \in \mathcal{P}_d^{(\mathrm{lc})}$, a $(1 - \alpha)$-confidence set for the distribution $P$ is given by
$$
	C_{\alpha,n}^{} = C_{\alpha,n}^{}(X_1,X_2,\ldots,X_n)
	\ := \ \Bigl\{ Q \in \mathcal{P}_d^{(\mathrm{lc})} :
		\bigl\| Q - \hat{P}_n \bigr\|_{\mathcal{H}} \le n^{-1/2} c_{\alpha,d} \Bigr\} .
$$
This entails simultaneous $(1 - \alpha)$-confidence sets for all integrals $\int \Pi(x) \, P(dx)$, where $\Pi: \Rd \to \RR$ is an arbitrary polynomial, namely,
$$
	C_{\alpha,n}^{(\Pi)} = C_{\alpha,n}^{(\Pi)}(X_1,X_2,\ldots,X_n)
	\ := \ \biggl\{ \int \Pi(x) \, Q(dx) : Q \in C_{\alpha,n}^{} \biggr\} .
$$
Since convergence with respect to $\|\cdot\|_{\mathcal{H}}$ implies weak convergence, Theorem~\ref{thm:aebe} implies the consistency of the confidence sets $C_{\alpha,n}^{(\Pi)}$, in the sense that
$$
	\sup_{t \in C_{\alpha,n}^{(\Pi)}} \Bigl| t - \int \Pi(x) \, P(dx) \Bigr|
	\to_p 0
	\quad\text{as} \ n \to \infty .
$$

Note that this construction proves existence of honest simultaneous confidence sets for arbitrary moments. But their explicit computation requires substantial additional work and is beyond the scope of the present paper.

\section{Various inequalities for $\mathcal{P}_d^{(\mathrm{lc})}$}
\label{sec: Various inequalities}

In this section we provide a few inequalities for log-concave distributions which are essential for the main result or are of independent interest. Let us first introduce some notation. The convex hull of a nonvoid set $S \subset \Rd$ is denoted by $\mathrm{conv}(S)$, the Lebesgue measure of a Borel set $S \subset \Rd$ by $|S|$.

\subsection{Inequalities for general dimension}
\label{subsec: Inequalities for general d}

\begin{lem}
\label{lem:andre}
Let $P \in \mathcal{P}_d^{(\mathrm{lc})}$ with density $f$. Let $x_0, x_1, \ldots, x_d$ be fixed points in $\Rd$ such that $\Delta := \conv\{x_0,x_1,\ldots,x_d\}$ has nonvoid interior. Then
$$
	\prod_{j=0}^d f(x_j)
	\ \le \ \Bigl( \frac{P(\Delta)}{|\Delta|} \Bigr)^{d+1} .
$$
Suppose that $x_1, x_2, \ldots, x_d \in \{f > 0\}$, and define $\tilde{f}(x_1,\ldots,x_d) := \Bigl( \prod_{i=1}^d f(x_i) \Bigr)^{1/d}$. Then
$$
	\frac{f(x_0)}{\tilde{f}(x_1,\ldots,x_d)}
	\ \le \ \Bigl(
		\frac{P(\Delta)}{\tilde{f}(x_1,\ldots,x_d) |\Delta|} \Bigr)^{d+1}
$$
If the right hand side is less than or equal to one, then
$$
	\frac{f(x_0)}{\tilde{f}(x_1,\ldots,x_d)}
	\ \le \ \exp \Bigl( 
		d - d \, \frac{\tilde{f}(x_1,\ldots,x_d) |\Delta|}{P(\Delta)} \Bigr) .
$$
\end{lem}

This lemma entails various upper bounds including a subexponential tail bound for log-concave densities.

\begin{lem}
\label{lem:dominic}
Let $x_0, x_1, \ldots, x_d \in \Rd$ and $\Delta$ as in Lemma~\ref{lem:andre}. Then for any $P \in \mathcal{P}_d^{(\mathrm{lc})}$ with density $f$ such that $x_0,x_1,\ldots,x_d \in \{f > 0\}$ and arbitrary $y \in \Delta$,
$$
	\min_{i=0,\ldots,d} f(x_i)
	\ \le \ f(y) \ \le \ \biggl( \frac{P(\Delta)}{|\Delta|} \biggr)^{d+1}
		\Bigl( \min_{i=0,\ldots,d} f(x_i) \Bigr)^{-d} .
$$
\end{lem}

\begin{lem}
\label{lem:lutz}
Let $x_0, x_1, \ldots, x_d \in \Rd$ as in Lemma~\ref{lem:andre}. Then there exists a constant $C = C(x_0, x_1, \ldots, x_d) > 0$ with the following property: For any $P \in \mathcal{P}_d^{(\mathrm{lc})}$ with density $f$ such that $x_0, x_1, \ldots, x_d \in \{f > 0\}$ and arbitrary $y \in \Rd$,
$$
	f(y) \ \le \ \max_{i=0,\ldots,d} f(x_i) \,
		H \Bigl( C \min_{i=0,\ldots,d} f(x_i) \, (1 + \|y\|^2)^{1/2} \Bigr) ,
$$
where
$$
	H(t) \ := \ \left\{\begin{array}{cl}
		t^{-(d+1)} & \text{for} \ t \in [0,1] , \\
		\exp(d - dt) & \text{for} \ t \ge 1 .
	\end{array}\right.
$$
\end{lem}

\begin{cor}
\label{cor:lutz}
For any $P \in \mathcal{P}_d^{(\mathrm{lc})}$ with density $f$ there exist constants $C_1 = C_1(P) > 0$ and $C_2 = C_2(P) > 0$ such that
$$
	f(x) \ \le \ C_1 \exp(- C_2 \|x\|)
	\quad\text{for all} \ x \in \Rd .
$$
\end{cor}

\subsection{Inequalities for dimension one}
\label{subsec: Inequalities for d=1}

In the special case $d = 1$ we denote the cumulative distribution function of $P$ with $F$. The hazard functions $f/F$ and $f/(1 - F)$ have the following properties:

\begin{lem}
\label{lem:dissandre1}
The function $f/F$ is non-increasing on $\{x : 0 < F(x) \le 1\}$, and the function $f/(1 - F)$ is non-decreasing on $\{x : 0 \le F(x) < 1\}$.\\[1ex]
Let $t_\ell := \inf \{f > 0\}$ and $t_u := \sup \{f > 0\}$. Then
\begin{eqnarray*}
	\lim_{t \downarrow t_\ell} \, \frac{f(t)}{F(t)} & = & \infty	\quad\text{if} \ t_\ell > - \infty , \\
	\lim_{t \uparrow t_u} \, \frac{f(t)}{1 - F(t)} & = & \infty	\quad\text{if} \ t_u < \infty .
\end{eqnarray*}
\end{lem}

\noindent
The monotonicity properties of the hazard functions $f/F$ and $f/(1-F)$ have been noted by An (1998) \nocite{An_1998} and Bagnoli and Bergstrom (2005) \nocite{Bagnoli_Bergstrom_2005}. For the reader's convenience a complete proof of Lemma~\ref{lem:dissandre1} will be given.

The next lemma provides an inequality for $f$ in terms of its first and second moments:

\begin{lem}
\label{lem:dissandre2}
Let $\mu$ and $\sigma$ be the mean and standard deviation, respectively, of the distribution $P$. Then for arbitrary $x_o \in \RR$,
$$
	f(x_o)^2 \ \le \ \frac{2 F(x_o)^3 + 2 (1 - F(x_o))^3}{(x_o - \mu)^2 + \sigma^2} .
$$
Equality holds if, and only if, $f$ is log-linear on both $(-\infty,x_o]$ and $[x_o,\infty)$.
\end{lem}

\section{Proofs}
\label{sec: Proofs}

\subsection{Proofs for Section~\ref{sec: Various inequalities}}

Our proof of Lemma~\ref{lem:andre} is based on a particular representation of Lebesgue measure on simplices: Let
$$
	\Delta_o \ := \ \bigl\{ u \in [0,1]^d : \sum_{i=1}^d u_i \le 1 \bigr\} .
$$
Then for any measurable function $h : \Delta_o \to [0,\infty)$,
$$
	\int_{\Delta_o} h(u) \, du
	\ = \ \frac{1}{d!} \, \EE \, h(B_1, B_2, \ldots, B_d) ,
$$
where $B_i := E_i \Big/ \sum_{j=0}^d E_j$ with independent, standard exponentially distributed random variables $E_0, E_1, \ldots, E_d$. This follows from general considerations about gamma and multivariate beta distributions, e.g.\ in Cule and D\"umbgen (2008). In particular, $|\Delta_o| = 1/d!$. Moreover, each variable $B_i$ is beta distributed with parameters $1$ and $d$, and $\EE(B_i) = 1/(d+1)$.

\paragraph{Proof of Lemma~\ref{lem:andre}.}
Any point $x \in \Delta$ may be written as
$$
	x(u) \ := \ x_0 + \sum_{i=1}^d u_i (x_i - x_0)
	\ = \ \sum_{i=0}^d u_i x_i
$$
for some $u \in \Delta_o$, where $u_0 := 1 - \sum_{i=1}^d u_i$. In particular,
$$
	\frac{|\Delta|}{|\Delta_o|}
	\ = \ \bigl| \det(x_1 - x_0, x_2 - x_0, \ldots, x_d - x_0) \bigr| .
$$
By concavity of $\varphi := \log f$,
$$
	\varphi(x(u)) \ \ge \ \sum_{i=0}^d u_i \varphi(x_i)
$$
for any $u = (u_i)_{i=1}^d \in \Delta_o$ and $u_0 = 1 - \sum_{i=1}^d u_i$. Hence
$$
	\frac{P(\Delta)}{|\Delta|}
	\ = \ \frac{1}{|\Delta_o|}
		\int_{\Delta_o} \exp \bigl( \varphi(x(u)) \bigr) \, du
	\ = \ \EE \exp \Bigl( \varphi \Bigl( \sum_{i=0}^d B_i x_i \Bigr) \Bigr)
	\ \ge \ \EE \exp \Bigl( \sum_{i=0}^d B_i \varphi(x_i) \Bigr) ,
$$
and by Jensen's inequality, the latter expected value is not less than
$$
	\exp \Bigl( \sum_{i=0}^d \EE(B_i) \varphi(x_i) \Bigr)
	\ = \ \exp \Bigl( \frac{1}{d+1} \sum_{i=0}^d \varphi(x_i) \Bigr)
	\ = \ \biggl( \prod_{i=0}^d f(x_i) \biggr)^{1/(d+1)} .
$$
This yields the first assertion of the lemma.

The inequality $\prod_{i=0}^d f(x_i) \le \bigl( P(\Delta)/|\Delta| \bigr)^{d+1}$ may be rewritten as
$$
	f(x_0) \tilde{f}(x_1,\ldots,x_d)^d
	\ \le \ \Bigl( \frac{P(\Delta)}{|\Delta|} \Bigr)^{d+1} ,
$$
and dividing both sides by $\tilde{f}(x_1,\ldots,x_d)^{d+1}$ yields the second assertion.

As to the third inequality, suppose that $f(x_0) \le \tilde{f}(x_1,\ldots,x_d)$, which is equivalent to $\varphi_0 := \varphi(x_0)$ being less than or equal to $\bar{\varphi} := \log \tilde{f}(x_1,\ldots,x_d) = d^{-1} \sum_{i=1}^d \varphi(x_i)$. Then
$$
	\frac{P(\Delta)}{|\Delta|}
	\ \ge \ \EE \exp \Bigl( \sum_{i=0}^d B_i \varphi(x_i) \Bigr)
	\ = \ \EE \exp \Bigl( B_0 \varphi_0
		+ (1 - B_0) \sum_{i=1}^d \tilde{B}_i \varphi(x_i) \Bigr) ,
$$
where $\tilde{B}_i := E_i \big/ \sum_{j=1}^d E_j$ for $1 \le i \le d$. It is well-known (e.g.\ Cule and D\"umbgen 2008\nocite{Cule_Duembgen_2008}) that $B_0$ and $\bigl( \tilde{B}_i \bigr)_{i=1}^d$ are stochastically independent, where $\EE \bigl( \tilde{B}_i \bigr) = 1/d$. Hence it follows from Jensen's inequality and $B_0 \sim \mathrm{Beta}(1,d)$ that
\begin{eqnarray*}
	\frac{P(\Delta)}{|\Delta|}
	& \ge & \EE \, \EE \biggl( \exp \Bigl( B_0 \varphi_0
		+ (1 - B_0) \sum_{i=1}^d \tilde{B}_i \varphi(x_i) \Bigr)
		\, \bigg| \, B_0 \biggr) \\
	& \ge & \EE \, \exp \biggl( \EE \Bigl( B_0 \varphi_0
			+ (1 - B_0) \sum_{i=1}^d \tilde{B}_i \varphi(x_i) \,\Big|\, B_0 \Bigr)
		\biggr) \\
	& = & \EE \, \exp \bigl( B_0 \varphi_0 + (1 - B_0) \bar{\varphi} \bigr) \\
	& = & \int_0^1 d (1 - t)^{d-1}
		\exp \bigl( t \varphi_0 + (1 - t) \bar{\varphi} \bigr) \, dt \\
	& = & \tilde{f}(x_1,\ldots,x_d)
		\int_0^1 d (1 - t)^{d-1}
			\exp \bigl( - t (\bar{\varphi} - \varphi_0) \bigr) \, dt \\
	& \ge & \tilde{f}(x_1,\ldots,x_d)
		\int_0^1 d (1 - t)^{d-1}
			\exp \bigl( \log(1 - t) (\bar{\varphi} - \varphi_0) \bigr) \, dt \\
	& = & \tilde{f}(x_1,\ldots,x_d)
		\int_0^1 d (1 - t)_{}^{\bar{\varphi} - \varphi_0 + d-1} \, dt \\
	& = & \tilde{f}(x_1,\ldots,x_d) \, \frac{d}{d + \bar{\varphi} - \varphi_0} .
\end{eqnarray*}
Thus $\bar{\varphi} - \varphi_0 \ge d \tilde{f}(x_1,\ldots,x_d) |\Delta| / P(\Delta) - d$, which is equivalent to
$$
	\frac{f(x_0)}{\tilde{f}(x_1,\ldots,x_d)}
	\ \le \ \exp \Bigl( d - d \, \frac{\tilde{f}(x_1,\ldots,x_d) |\Delta|}{P(\Delta)}
		\Bigr) .
	\eqno{\Box}
$$

We first prove Lemma~\ref{lem:lutz} because this provides a tool for the proof of Lemma~\ref{lem:dominic} as well.

\paragraph{Proof of Lemma~\ref{lem:lutz}.}
At first we investigate how the size of $\Delta$ changes if we replace one of its vertices with another point. Note that for any fixed index $j \in \{0,1,\ldots,d\}$,
$$
	\bigl| \det(x_i - x_j : i \ne j) \bigr|
	\ = \ |\det(X)|
	\quad\text{with}\quad
	X \ := \ \biggl(\!\!\begin{array}{cccc}
			x_0 & x_1 & \ldots & x_d \\
			1   & 1   & \ldots & 1
		\end{array}\!\!\biggr) .
$$
Moreover, any point $y \in \Rd$ has a unique representation $y = \sum_{i=0}^d \lambda_i x_i$ with scalars $\lambda_0$, $\lambda_1$, \ldots, $\lambda_d$ summing to one. Namely,
$$
	(\lambda_i)_{i=0}^d
	\ = \ X^{-1} \biggl(\!\!\begin{array}{c} y \\ 1 \end{array}\!\!\biggr) .
$$
Hence the set $\Delta_j(y) := \conv \bigl( \{x_i : i \ne j\} \cup \{y\} \bigr)$ has Lebesgue measure
\begin{eqnarray*}
	|\Delta_j(y)|
	& = & \frac{1}{d!} \, \biggl| \det
		\biggl(\!\!\begin{array}{ccccccc}
			x_0 & \ldots & x_{j-1} & y & x_{j+1} & \ldots & x_d \\
			1   & \ldots & 1       & 1 & 1       & \ldots & 1
		\end{array}\!\!\biggr) \biggr| \\
	& = & \frac{1}{d!} \, \biggl| \sum_{i=0}^d \lambda_i \det
		\biggl(\!\!\begin{array}{ccccccc}
			x_0 & \ldots & x_{j-1} & x_i & x_{j+1} & \ldots & x_d \\
			1   & \ldots & 1       & 1   & 1       & \ldots & 1
		\end{array}\!\!\biggr) \biggr| \\
	& = & \frac{1}{d!} \, |\lambda_j| |\det(X)| \\
	& = & |\lambda_j| |\Delta| .
\end{eqnarray*}
Consequently,
\begin{eqnarray*}
	\max_{j=0,1,\ldots,d} |\Delta_j(y)|
	& = & |\Delta| \, \max_{j=0,1,\ldots,d} |\lambda_j| \\
	& = & |\Delta| \, \biggl\|
			X^{-1} \biggl(\!\!\begin{array}{c} y \\ 1 \end{array}\!\!\biggr)
		\biggr\|_\infty \\
	& \ge & |\Delta| (d+1)_{}^{-1/2} \biggl\|
			X^{-1} \biggl(\!\!\begin{array}{c} y \\ 1 \end{array}\!\!\biggr)
		\biggr\| \\
	& \ge & |\Delta| (d+1)_{}^{-1/2} \sigma_{\rm max}(X)^{-1} (\|y\|^2 + 1)^{1/2} ,
\end{eqnarray*}
where $\sigma_{\rm max}(X) > 0$ is the largest singular value of $X$.

Now we consider any log-concave probability density $f$. Let $f_{\rm min}$ and $f_{\rm max}$ denote the minimum and maximum, respectively, of $\{f(x_i) : i=0,\ldots,d\}$, where $f_{\rm min}$ is assumed to be greater than zero. Applying Lemma~\ref{lem:andre} to $\Delta_j(y)$ in place of $\Delta$ with suitably chosen index $j$, we may conclude that
$$
	f(y) \ \le \ f_{\rm max} \bigl( C f_{\rm min} (\|y\|^2 + 1)^{1/2} \bigr)^{-(d+1)} ,
$$
where $C = C(x_0,\ldots, x_d) := |\Delta| (d+1)_{}^{-1/2} \sigma_{\rm max}(X)^{-1}$. Moreover, in case of $C f_{\rm min} (\|y\|^2 + 1)^{1/2} \ge 1$,
$$
	f(y) \ \le \ f_{\rm max} \exp \bigl( d - d C f_{\rm min} (\|y\|^2 + 1)^{1/2} \bigr) .
	\eqno{\Box}
$$

\paragraph{Proof of Lemma~\ref{lem:dominic}.}
Let $y \in \Delta$, i.e.\ $y = \sum_{i=0}^d \lambda_i x_i$ with a unique vector $\lambda = (\lambda_i)_{i=0}^d$ in $[0,1]^{d+1}$ whose components sum to one. With $\Delta_j(y)$ as in the proof of Lemma~\ref{lem:lutz}, elementary calculations reveal that
$$
	\Delta \ = \ \bigcup_{j \in J} \Delta_j(y) ,
$$
where $J := \{j : \lambda_j > 0\}$. Moreover, all these simplices $\Delta_j(y)$, $j \in J$, have nonvoid interior, and $|\Delta_j(y) \cap \Delta_k(y)| = 0$ for different $j,k \in J$. Consequently it follows from Lemma~\ref{lem:andre} that
\begin{eqnarray*}
	\frac{P(\Delta)}{|\Delta|}
	& = & \sum_{j \in J}
		\frac{|\Delta_j(y)|}{|\Delta|} \cdot \frac{P(\Delta_j(y))}{|\Delta_j(y)|} \\
	& \ge & \sum_{j \in J}
		\frac{|\Delta_j(y)|}{|\Delta|} \cdot \Bigl( f(y) \prod_{i \ne j} f(x_i) \Bigr)^{1/(d+1)} \\
	& \ge & \sum_{j \in J}
		\frac{|\Delta_j(y)|}{|\Delta|} \cdot f(y)^{1/(d+1)}
			\Bigl( \min_{i=0,\ldots,d} f(x_i) \Bigr)^{d/(d+1)} \\
	& = & f(y)^{1/(d+1)} \Bigl( \min_{i=0,\ldots,d} f(x_i) \Bigr)^{d/(d+1)} .
\end{eqnarray*}
This entails the asserted upper bound for $f(y)$. The lower bound follows from the elementary fact that any concave function on the simplex $\Delta$ attains its minimal value in one of the vertices $x_0, x_1, \ldots, x_d$.	\hfill	$\Box$

\paragraph{Proof of Lemma~\ref{lem:dissandre1}.}
We only prove the assertions about $f/(1 - F)$. Considering the distribution function $\tilde{F}(x) := 1 - F(- x)$ with log-concave density $\tilde{f}(x) = f(- x)$ then yields the corresponding properties of $f/F$.

Note that $\{F < 1\} = (-\infty, t_u)$. On $\{f = 0\} \cap (-\infty,t_u)$, the function $f/(1 - F)$ is equal to zero. For $t \in \{f > 0\} \cap (-\infty,t_u)$,
$$
	\frac{f(t)}{1 - F(t)}
	\ = \ \Bigl( \int_0^\infty \exp \bigl( \varphi(t+x) - \varphi(t) \bigr) \, dx \Bigr)^{-1}
$$
is non-decreasing in $t$, because $t \mapsto \varphi(t+x) - \varphi(t)$ is non-increasing in $t \in \{f > 0\}$ for any fixed $x > 0$, due to concavity of $\varphi$.

In case of $t_u < \infty$, fix any point $s \in (t_\ell, t_u)$. Then for $s \le t < t_u$,
\begin{eqnarray*}
	\frac{f(t)}{1 - F(t)}
	& = & \Bigl( \int_t^{t_u} \exp \bigl( \varphi(x) - \varphi(t) \bigr) \, dx \Bigr)^{-1} \\
	& \ge & \Bigl( \int_t^{t_u} \exp \bigl( \varphi'(s\,+) (x - t) \bigr) \, dx \Bigr)^{-1} \\
	& \ge & \Bigl( \exp \bigl( \min(\varphi'(s\,+), 0) (t_u - t) \bigr) (t_u - t) \Bigr)^{-1} \\
	& \to & \infty	\quad\text{as} \ t \uparrow t_u .
\end{eqnarray*}\\[-7ex]
\strut \hfill	$\Box$

\paragraph{Proof of Lemma~\ref{lem:dissandre2}.}
The asserted upper bound for $f(t_o)$ is strictly positive and continuous in $t_o$. Hence it suffices to consider a point $t_o$ with $0 < F(t_o) < 1$. Since $(x_o - \mu)^2 + \sigma^2$ equals $\int (x - x_o)^2 f(x) \, dx$, we try to bound the latter integral from above. To this end, let $g$ be a piecewise loglinear probability density, namely,
$$
	g(x) \ := \ \begin{cases}
		f(x_o) \exp( - a |x - x_o|) & \text{if} \ x \le x_o , \\
		f(x_o) \exp( - b |x - x_o|) & \text{if} \ x \ge x_o ,
	\end{cases}
$$
with $a := f(x_o) / F(x_o)$ and $b := f(x_o) / (1 - F(x_o))$, so that
$$
	\int_{-\infty}^{x_o} (g - f)(x) \, dx \ = \ \int_{x_o}^{\infty} (g - f)(x) \, dx \ = \ 0 .
$$
By concavity of $\log f$, there are real numbers $r < x_o < s$ such that $f \ge g$ on $(r,s)$ and $f \le g$ on $\RR \setminus [r,s]$. Consequently,
\begin{eqnarray*}
	\int (x - x_o)^2 (f - g)(x) \, dx
	& = & \int_{-\infty}^{x_o} \underbrace{\bigl[ (x - x_o)^2 - (r - x_o)^2 \bigr] (f - g)(x)}_{\le \ 0} \, dx \\
	&& + \ \int_{x_o}^{\infty} \underbrace{\bigl[ (x - x_o)^2 - (s - x_o)^2 \bigr] (f - g)(x)}_{\le \ 0} \, dx \\
	& \le & 0 ,
\end{eqnarray*}
with equality if, and only if, $f = g$. Now the assertion follows from
\begin{eqnarray*}
	\int (x - x_o)^2 g(x) \, dx
	& = & f(x_o) \Bigl( \int_0^\infty t^2 \exp(- at) \, dt + \int_0^\infty t^2 \exp(- bt) \, dt \Bigr) \\
	& = & \frac{2 F(x_o)^3 + 2 (1 - F(x_o))^3}{f(x_o)^2} .
\end{eqnarray*}\\[-7ex]
\strut \hfill	$\Box$

\subsection{Proof of the main results}

Note first that $\{f > 0\}$ is a convex set with nonvoid interior. For notational convenience we may and will assume that
$$
	0 \ \in \ \interior \{f > 0\} .
$$
For if $x_o$ is any fixed interior point of $\{f > 0\}$ we could just shift the coordinate system and consider the densities $\tilde{f} := f(x_o + \cdot)$ and $\tilde{f}_n := f_n(x_o + \cdot)$ in place of $f$ and $f_n$, respectively. Note also that $A(x_o + x) - A(x) \in \bigl[ - A(- x_o), A(x_o) \bigr]$, due to subadditivity of $A$.

In our proof of Theorem~\ref{thm:dabigone}, Part~(i), we utilize two simple inequalities for log-concave densities:

\begin{lem}
\label{lem:dasmallone1}
Let $x_0, x_1, \ldots, x_d \in \Rd$ such that $\Delta := \conv\{x_0,x_1,\ldots,x_d\}$ has nonvoid interior. For $j=0,1,\ldots,d$ define the ``corner simplex''
$$
	\Delta_j \ := \ \bigl\{ 2 x_j - x : x \in \Delta \} ,
$$
i.e.\ the reflection of $\Delta$ at the point $x_j$. Let $P \in \mathcal{P}_d^{(\mathrm{lc})}$ with density $f = \exp \compo \varphi$. If $P(\Delta_j) > 0$ for all $j=0,1,\ldots,d$, then $\Delta \subset \interior \{f > 0\}$, and
\begin{eqnarray*}
	\min_{j=0,1,\ldots,d} \log \frac{P(\Delta_j)}{|\Delta|}
	& \le & \min_{x \in \Delta} \varphi(x)
		\ \le \ \log \frac{P(\Delta)}{|\Delta|} \\
	& \le & \max_{x \in \Delta} \varphi(x)
		\ \le \ (d+1) \log \frac{P(\Delta)}{|\Delta|}
			- d \min_{j=0,1,\ldots,d} \log \frac{P(\Delta_j)}{|\Delta|} .
\end{eqnarray*}
\end{lem}

\begin{figure}[h]
\centerline{\includegraphics[width=0.7\textwidth]{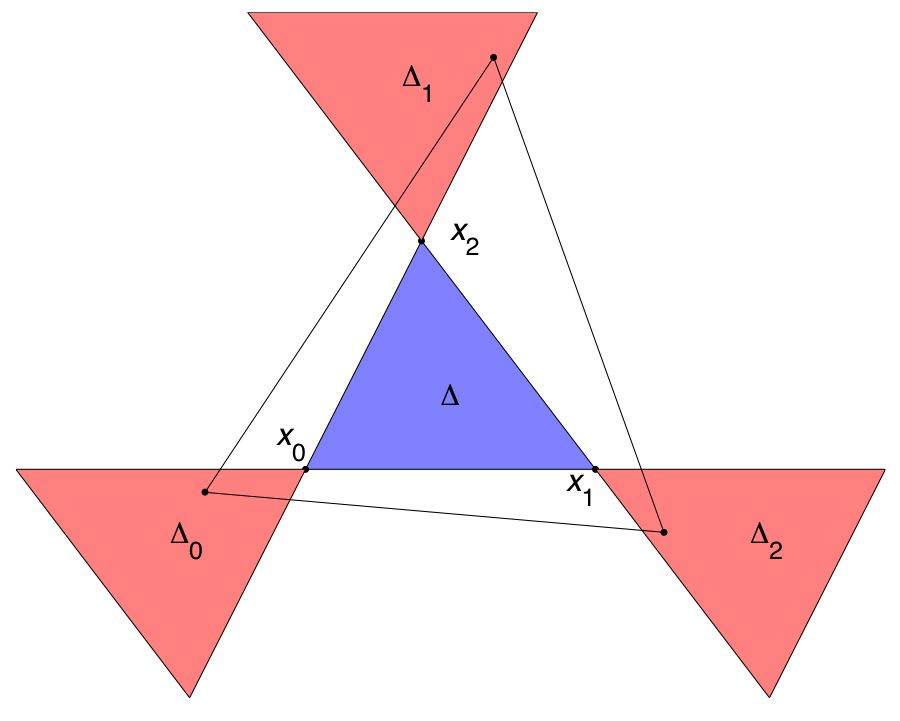}}
\caption{A simplex $\Delta$ and its corner simplices $\Delta_j$.}
\label{fig:CornerSimplices}
\end{figure}

\noindent
Figure~\ref{fig:CornerSimplices} illustrates the definition of the corner simplices and a key statement in the proof of Lemma~\ref{lem:dasmallone1}.

\begin{lem}
\label{lem:dasmallone2}
Suppose that $B(0,\delta) \subset \{f > 0\}$ for some $\delta > 0$. For $t \in (0,1)$ define $\delta_t := (1 - t) \delta /(1 + t)$. Then for any $y \in \Rd$,
$$
	\sup_{x \in B(y, \delta_t)} f(x)
	\ \le \ \Bigl( \inf_{v \in B(0,\delta)} f(v) \Bigr)^{1 - 1/t}
	\Bigl( \frac{P(B(ty, \delta_t)}{|B(ty,\delta_t)|} \Bigr)^{1/t} .
$$
\end{lem}

\begin{figure}[h]
\centerline{\includegraphics[width=0.6\textwidth]{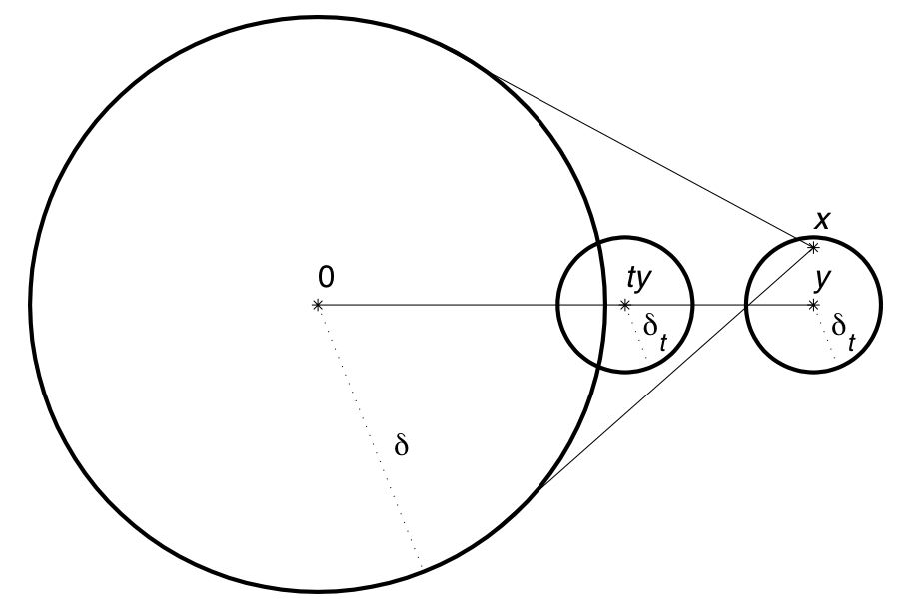}}
\caption{The three closed balls in Lemma~\ref{lem:dasmallone2}.}
\label{fig:ThreeBalls}
\end{figure}

\noindent
This lemma involves three closed balls $B(0,\delta)$, $B(ty, \delta_t)$ and $B(y, \delta_t)$; see Figure~\ref{fig:ThreeBalls} for an illustration of these and the key argument of the proof.

\paragraph{Proof of Lemma~\ref{lem:dasmallone1}.}
Suppose that all corner simplices satisfy $P(\Delta_j) > 0$. Then for $j=0,1,\ldots,d$ there exists an interior point $z_j$ of $\Delta_j$ with $f(z_j) > 0$, that means, $z_j = 2x_j - \sum_{i=0}^d \lambda_{ij} x_i$ with positive numbers $\lambda_{ij}$ such that $\sum_{i=0}^d \lambda_{ij} = 1$. With the matrices
$$
	X \ := \ \begin{pmatrix} x_0 & x_1 & \ldots & x_d \\ 1 & 1 & \ldots & 1 \end{pmatrix} ,
	\quad
	Z \ := \ \begin{pmatrix} z_0 & z_1 & \ldots & z_d \\ 1 & 1 & \ldots & 1 \end{pmatrix}
	\quad\text{and}\quad
	\Lambda \ := \ \begin{pmatrix}
		\lambda_{00} & \ldots & \lambda_{0d} \\
		\vdots & & \vdots \\
		\lambda_{d0} & \ldots & \lambda_{dd}
	\end{pmatrix}
$$
in $\mathbb{R}^{(d+1)\times(d+1)}$ we may write
$$
	Z \ = \ X (2 I - \Lambda) .
$$
But the matrix $2I - \Lambda$ is nonsingular with inverse
$$
	M \ := \ (2I - \Lambda)^{-1}
	\ = \ 2^{-1} (I - 2^{-1} \Lambda)^{-1}
	\ = \ \sum_{\ell=0}^\infty 2^{-(\ell+1)} \Lambda^\ell .
$$
The latter power series converges, because $\Lambda^\ell$ has positive components for all $\ell \ge 1$, and via induction on $\ell \ge 0$ one can show that all columns of $\Lambda^\ell$ sum to one. Consequently, $X = Z M$, i.e.\ for each index $j$, the point $x_j$ may be written as $\sum_{i=0}^d \mu_{ij} z_i$ with positive numbers $\mu_{ij}$ such that $\sum_{i=0}^d \mu_{ij} = 1$. This entails that $\Delta$ is a subset of $\interior\conv\{z_0,z_1,\ldots,z_d\} \subset \interior \{f > 0\}$; see also Figure~\ref{fig:CornerSimplices}.

Since $\min_{x \in \Delta} f(x) \le P(\Delta)/|\Delta| \le \max_{x \in \Delta} f(x)$, the inequalities
$$
	\min_{x \in \Delta} \varphi(x) \ \le \ \log \frac{P(\Delta)}{|\Delta|}
	\ \le \ \max_{x \in \Delta} \varphi(x)
$$
are obvious. By concavity of $\varphi$, its minimum over $\Delta$ equals $\varphi(x_{j_o})$ for some index $j_o \in \{0,1,\ldots,d\}$. But then for arbitrary $x \in \Delta$ and $y := 2x_{j_o} - x \in \Delta_{j_o}$, it follows from $x_{j_o} = 2^{-1}(x + y)$ and concavity of $\varphi$ that
$$
	\varphi(x_{j_o}) \ \ge \ \frac{\varphi(x) + \varphi(y)}{2}
	\ \ge \ \frac{\varphi(x_{j_o}) + \varphi(y)}{2} ,
$$
so that $\varphi \le \varphi(x_{j_o})$ on $\Delta_{j_o}$. Hence
$$
	\min_{x \in \Delta} \varphi(x) \ = \ \varphi(x_{j_o})
	\ \ge \ \log \frac{P(\Delta_{j_o})}{|\Delta|} .
$$
Finally, Lemma~\ref{lem:dominic} entails that
\begin{eqnarray*}
	\max_{x \in \Delta} \varphi(x)
	& \le & (d+1) \log \frac{P(\Delta)}{|\Delta|}
		- d \min_{j=0,1,\ldots,d} \varphi(x_j) \\
	& \le & (d+1) \log \frac{P(\Delta)}{|\Delta|}
		- d \min_{j=0,1,\ldots,d} \log \frac{P(\Delta_j)}{|\Delta|} .
\end{eqnarray*}\\[-7ex]
\strut	\hfill	$\Box$

\paragraph{Proof of Lemma~\ref{lem:dasmallone2}.}
The main point is to show that for any point $x \in B(y, \delta_t)$,
$$
	B(ty, \delta_t) \ \subset \ (1 - t) B(0, \delta) + t x ,
$$
i.e.\ any point $w \in B(ty, \delta_t)$ may be written as $(1 - t) v + t x$ for a suitable $v \in B(0, \delta)$; see also Figure~\ref{fig:ThreeBalls}. But note that the equation $(1 - t) v + t x = w$ is equivalent to $v = (1 - t)^{-1}(w - tx)$. This vector $v$ belongs indeed to $B(0,\delta)$, because
$$
	\|v\| \ = \ (1 - t)^{-1} \|w - tx\|
	\ = \ (1 - t)^{-1} \bigl\| w - ty + t(y - x) \bigr\|
	\ \le \ (1 - t)^{-1} (\delta_t + t \delta_t)
	\ = \ \delta
$$
by definition of $\delta_t$.

This consideration shows that for any point $x \in B(y, \delta_t)$ and any point $w \in B(ty, \delta_t)$,
$$
	f(w) \ \ge \ f(v)^{1-t} f(x)^t
	\ \ge \ J_0^{1 - t} f(x)^t
$$
with $v = (1 - t)^{-1} (w - tx) \in B(0, \delta)$ and $J_0 := \inf_{v \in B(0,\delta)} f(v)$. Averaging this inequality with respect to $w \in B(ty, \delta_t)$ yields
$$
	\frac{P(B(ty,\delta_t))}{|B(ty,\delta_t)|}
	\ \ge \ J_0^{1 - t} f(x)^t .
$$
Since $x \in B(y,\delta_t)$ is arbitrary, this entails the assertion of Lemma~\ref{lem:dasmallone2}.	\hfill	$\Box$

\paragraph{Proof of Theorem~\ref{thm:dabigone}, Part~(i).}
Our proof is split into three steps.

{\sl
\paragraph{Step 1:} The sequence $(f_n)_n$ converges to $f$ uniformly on any compact subset of $\interior \{f > 0\}$.
} 

\noindent
By compactness, this claim is a consequence of the following statement: For any interior point $y$ of $\{f > 0\}$ and any $\eta > 0$ there exists a neighborhood $\Delta(y,\eta)$ of $y$ such that
$$
	\limsup_{n \to \infty} \ \sup_{x \in \Delta(y,\eta)} \Bigl| \frac{f_n(x)}{f(x)} - 1 \Bigr|
	\ \le \ \eta .
$$
To prove the latter statement, fix any number $\epsilon \in (0,1)$. Since $f$ is continuous on $\interior \{f > 0\}$, there exists a simplex $\Delta = \conv\{x_0,x_1,\ldots,x_d\}$ such that $y \in \interior \Delta$ and
$$
	f \ \in \ \bigl[ (1 - \epsilon) f(y), (1 + \epsilon) f(y) \bigr]
	\quad\text{on}\quad \Delta \cup \Delta_0 \cup \Delta_1 \cup \cdots \cup \Delta_d
$$
with the corner simplices $\Delta_j$ defined as in Lemma~\ref{lem:dasmallone1}. Since the boundary of any simplex $\tilde{\Delta}$ is contained in the union of $d+1$ hyperplanes, it satisfies $P(\partial \tilde{\Delta}) = 0$, so that weak convergence of $(P_n)_n$ to $P$ implies that
$$
	\lim_{n \to \infty} P_n(\tilde{\Delta}) \ = \ P(\tilde{\Delta}) .
$$
Therefore it follows from Lemma~\ref{lem:dasmallone1} that
\begin{eqnarray*}
	\liminf_{n \to \infty} \ \inf_{x \in \Delta} \frac{f_n(x)}{f(x)}
	& \ge & \liminf_{n \to \infty} \ \frac{1}{(1+\epsilon) f(y)}
		\inf_{x \in \Delta} f_n(x) \\
	& \ge & \liminf_{n \to \infty} \ \frac{1}{(1+\epsilon) f(y)}
		\min_{j=0,1,\ldots,d} \frac{P_n(\Delta_j)}{|\Delta|} \\
	& = & \frac{1}{(1+\epsilon) f(y)}
		\min_{j=0,1,\ldots,d} \frac{P(\Delta_j)}{|\Delta|}
	\ \ge \ \frac{1-\epsilon}{1+\epsilon}
\end{eqnarray*}
and
\begin{eqnarray*}
	\limsup_{n \to \infty} \ \sup_{x \in \Delta} \frac{f_n(x)}{f(x)}
	& \le & \limsup_{n \to \infty} \ \frac{1}{(1-\epsilon) f(y)}
		\sup_{x \in \Delta} f_n(x) \\
	& \le & \frac{1}{(1-\epsilon) f(y)} \Bigl( \frac{P(\Delta)}{|\Delta|} \Bigr)^{d+1}
		\Bigl( \min_{j=0,1,\ldots,d} \frac{P(\Delta_j)}{|\Delta|} \Bigr)^{-d}
	\ \le \ \Bigl( \frac{1+\epsilon}{1-\epsilon} \Bigr)^{d+1} .
\end{eqnarray*}
For $\epsilon$ sufficiently small, both $(1 - \epsilon)/(1 + \epsilon) \ge 1 - \eta$ and $\bigl( (1+\epsilon)/(1 - \epsilon) \bigr)^{d+1} \le 1 + \eta$, which proves the assertion of step~1.

{\sl
\paragraph{Step 2:} If $f$ is continuous at $y \in \Rd$ with $f(y) = 0$, then for any $\eta > 0$ there exists a number $\delta(y,\eta) > 0$ such that
$$
	\limsup_{n \to \infty} \ \sup_{x \in B(y,\delta(y,\eta))} f_n(x) \ \le \ \eta \, . 
$$
} 

\noindent
For this step we employ Lemma~\ref{lem:dasmallone2}. Let $\delta_0 > 0$ such that $B(0,\delta_0)$ is contained in $\interior \{f>0\}$. Furthermore, let $J_0 > 0$ be the minimum of $f$ over $B(0,\delta_0)$. Then step~1 entails that
$$
	\liminf_{n \to \infty} \ \inf_{x \in B(0,\delta_0)} f_n(x) \ \ge \ J_0 .
$$
Moreover, for any $t \in (0,1)$ and $\delta_t := (1 - t) \delta_0 /(1 + t)$,
\begin{eqnarray*}
	\limsup_{n \to \infty} \ \sup_{x \in B(y,\delta_t)} \, f_n(x)
	& \le & J_0^{1 - 1/t} \
		\limsup_{n \to \infty}
			\Bigl( \frac{P_n(B(ty,\delta_t))}{|B(y,\delta_t)|} \Bigr)^{1/t} \\
	& \le & J_0^{1 - 1/t} \Bigl( \frac{P(B(ty,\delta_t))}{|B(y,\delta_t)|} \Bigr)^{1/t} \\
	& \le & J_0^{1 - 1/t} \Bigl( \sup_{x \in B(ty,\delta_t)} f(x) \Bigr)^{1/t} .
\end{eqnarray*}
But the latter bound tends to zero as $t \uparrow 1$.

{\sl
\paragraph{Final step:} $(f_n)_n$ converges to $f$ uniformly on any closed set of continuity points of $f$.
} 

\noindent
Let $S$ be such a closed set. Then Steps~1 and 2 entail that
$$
	\lim_{n \to \infty} \ \sup_{x \in S \cap B(0, \rho)} \bigl| f_n(x) - f(x) \bigr|
	\ = \ 0
$$
for any fixed $\rho \ge 0$, because $S \cap B(0,\rho)$ is compact, and any point $y \in S \setminus \interior \{f > 0\}$ satisfies $f(y) = 0$.

On the other hand, let $\Delta$ be a nondegenerate simplex with corners $x_0,x_1,\ldots,x_d \in \interior \{f > 0\}$. Step~1 also implies that $\lim_{n \to \infty} f_n(x_i) = f(x_i)$ for $i = 0,1,\ldots,d$, so that Lemma~\ref{lem:lutz} entails that \begin{equation}
\label{eq: tail bound}
	\limsup_{n \to \infty} \ \sup_{x \,:\, \|x\| \ge \rho} \max \bigl\{ f_n(x), f(x) \bigr\}
	\ \le \ \max_{i=0,\ldots,d} f(x_i)
		H \Bigl( C \min_{i=0,\ldots,d} f(x_i) (1 + \rho^2)^{1/2} \Bigr)
\end{equation}
for any $\rho \ge 0$ with a constant $C = C(x_0,\ldots,x_d) > 0$. Since this bound tends to zero as $\rho \to \infty$, the assertion of Theorem~\ref{thm:dabigone}, Part~(i) follows.	\hfill	$\Box$

Our proof of Theorem~\ref{thm:dabigone}, Part~(ii), is based on Part~(i) and an elementary result about convex sets:

\begin{lem}
\label{lem:dasmallone3}
Let $\mathcal{C}$ be a convex subset of $\Rd$ containing $B(0,\delta)$ for some $\delta > 0$. If $y \in \mathcal{C}$, then
$$
	B(ty, (1 - t)\delta) \ \subset \ \mathcal{C}
	\quad\text{for all} \ t \in [0,1] .
$$
If $y \in \Rd \setminus \mathcal{C}$, then
$$
	B(\lambda y, (\lambda - 1)\delta) \ \subset \ \Rd \setminus \mathcal{C}
	\quad\text{for all} \ \lambda \ge 1 .
$$
\end{lem}

\noindent
One consequence of this lemma is the well-known fact that the boundary of the convex set $\{f > 0\}$ has Lebesgue measure zero. Namely, for any unit vector $u \in \Rd$ there exists at most one number $r > 0$ such that $ru \in \partial \{f > 0\}$. Lemma~\ref{lem:dasmallone3} is needed to obtain a refinement of this fact.

\paragraph{Proof of Lemma~\ref{lem:dasmallone3}.}
By convexity of $\mathcal{C}$ and $B(0,\delta) \subset \mathcal{C}$, it follows from $y \in \mathcal{C}$ that
$$
	\mathcal{C} \ \supset \ \bigl\{ (1 - t)v + ty : v \in B(0,\delta) \bigr\}
	\ = \ B(ty, (1 - t)\delta)
$$
for any $t \in [0,1]$. In case of $y \not\in \mathcal{C}$, for $\lambda \ge 1$ and arbitrary $x \in B(\lambda y, (\lambda-1)\delta)$ we write $x = \lambda y + (\lambda - 1) v$ with $v \in B(0,\delta)$. But then
$$
	y \ = \ (1 - \lambda^{-1}) (-v) + \lambda^{-1} x .
$$
Hence $y \not\in \mathcal{C}$ is a convex combination of a point in $B(0,\delta) \subset \mathcal{C}$ and $x$, so that $x \not\in \mathcal{C}$, too. 	\hfill	$\Box$

\paragraph{Proof of Theorem~\ref{thm:dabigone}, Part~(ii).}
It follows from (\ref{eq: tail bound}) in the proof of Part~(i) with $\rho = 0$ that
$$
	\limsup_{n \to \infty} \ \sup_{x \in \Rd} \, f_n(x) \ < \ \infty .
$$
Since $(f_n)_n$ converges to $f$ pointwise on $\Rd \setminus \partial \{f > 0\}$, and since $\partial \{f > 0\}$ has Lebesgue measure zero, dominated convergence yields
\begin{eqnarray*}
	\lefteqn{ \limsup_{n \to \infty} \int_{\Rd} \exp(A(x)) \bigl| f_n(x) - f(x) \bigr| \, dx } \\
	& = & \limsup_{n \to \infty}
		\int_{\Rd \setminus B(0,\gamma)} \exp(A(x)) \bigl| f_n(x) - f(x) \bigr| \, dx \\
	& \le & \limsup_{n \to \infty}
		\int_{\Rd \setminus B(0,\gamma)} \exp(A(x)) \max \bigl( f_n(x),f(x) \bigr) \, dx
\end{eqnarray*}
for any fixed $\gamma > 0$.

It follows from Assumption~(\ref{eq:dacondition}) that for a suitable $\rho > 0$,
$$
	A(x) + \varphi(x) - \varphi(0) \ \le \ - 1
	\quad\text{whenever} \ \|x\| \ge \rho .
$$
Utilizing sublinearity of $A$ and concavity of $\varphi$, we may deduce that for $x \in \Rd$ with $\|x\| \ge \rho$ even
\begin{eqnarray*}
	A(x) + \varphi(x)
	& = & \varphi(0) + A(x) + \|x\| \frac{\varphi(\|x\| u) - \varphi(0)}{\|x\|} \\
	& \le & \varphi(0) + A(x) + \|x\| \frac{\varphi(\rho u) - \varphi(0)}{\rho} \\
	& = & \varphi(0) + \rho^{-1} \|x\| \bigl( A(\rho u) + \varphi(\rho u) - \varphi(0) \bigr) \\
	& \le & \varphi(0) - \rho^{-1} \|x\| ,
\end{eqnarray*}
where $u := \|x\|^{-1} x$.
In particular, $\int_{\Rd} \exp(A(x)) f(x) \, dx$ is finite. Now let $\delta > 0$ such that $B(0,\delta) \subset \{f > 0\}$. It follows from Lemma~\ref{lem:dasmallone3} that for any unit vector $u \in \Rd$, either $2\rho u \in \{f > 0\}$ and $B(\rho u, \delta/2) \subset \{f > 0\}$, or $2\rho u \in \{f = 0\}$ and $B(3\rho u, \delta/2) \subset \{f = 0\}$. Hence
$$
	K \ := \ \{0\} \cup \Bigl\{ x \in \Rd : \|x\| \in \{\rho, 3\rho\},
		\inf_{y \in \partial \{f > 0\}} \|x - y\| \ge \delta/2 \Bigr\}
$$
defines a compact subset of $\Rd \setminus \partial \{f > 0\}$ such that
$$
	K \cap \{\rho u, 3\rho u\} \ \neq \ \emptyset
	\quad\text{for any unit vector} \ u \in \Rd .
$$
According to Part~(i), $(f_n)_n$ converges to $f$ uniformly on $K$. Thus for fixed numbers $\epsilon' > 0$, $\epsilon'' \in (0,\rho^{-1})$ and sufficiently large $n$, the log-densities $\varphi_n := \log f_n$ satisfy the following inequalities:
\begin{eqnarray*}
	A(ru) + \varphi_n(ru)
	& = & \varphi_n(0)
		+ r \Bigl( A(u) + \frac{\varphi_n(ru) - \varphi_n(0)}{r} \Bigr) \\
	& \le & \varphi_n(0)
		+ r \Bigl( A(u) + \min_{s = \rho, 3\rho} \frac{\varphi_n(su) - \varphi_n(0)}{s} \Bigr) \\
	& \le & \varphi(0) + \epsilon' - \epsilon'' r
\end{eqnarray*}
for all unit vectors $u \in \Rd$ and $r \ge 3\rho$. Hence for $\gamma \ge 3\rho$,
\begin{eqnarray*}
	\lefteqn{ \limsup_{n \to \infty}
		\int_{\Rd \setminus B(0,\gamma)} \exp(A(x)) \max \bigl( f_n(x),f(x) \bigr) \, dx } \\
	& \le & f(0) \int_{\Rd \setminus B(0,\gamma)} \exp \bigl( \epsilon' - \epsilon'' \|x\|) \, dx \\
	& = & \mathrm{const}(d) f(0) \int_\gamma^\infty r^{d-1} \exp(\epsilon' - \epsilon'' r) \, dr \\
	& \to & 0	\quad\text{as} \ \gamma \to \infty .
\end{eqnarray*}\\[-7ex]
\strut	\hfill	$\Box$

\paragraph{Proof of Proposition~\ref{prop:godot}.}
It follows from convexity of $\exp(\cdot)$ that $\Theta(P)$ is a convex subset of $\Rd$, and obviously it contains $0$. Now we verify it to be open. For any fixed $\theta \in \Theta(P)$ we define a new probability density
$$
	\tilde{f}(x) \ := \ C^{-1} \exp(\theta^\top x) f(x)
	\ = \ \exp \bigl( \theta^\top x + \varphi(x) - \log C \bigr)
$$
with $C := \int_{\Rd} \exp(\theta^\top x) f(x) \, dx$. Obviously, $\tilde{f}$ is log-concave, too. Thus, by Corollary~\ref{cor:lutz}, there exist constants $C_1, C_2 > 0$ such that $\tilde{f}(x) \le C_1 \exp( - C_2 \|x\|)$ for all $x \in \Rd$. In particular,
$$
	\infty \ > \ C \int_{\Rd} \exp(\delta^\top x) \tilde{f}(x) \, dx
	\ = \ \int_{\Rd} \exp \bigl( (\theta + \delta)^\top x) f(x) \, dx
$$
for all $\delta \in \Rd$ with $\|\delta\| < C_2$. This shows that $\Theta(P)$ is open.

Finally, let $\theta \in \Theta(P)$ and $\epsilon > 0$ such that $B(\theta,\epsilon) \subset \Theta(P)$. With the previous arguments one can show that for each unit vector $u \in \Rd$ there exist constants $D(u) \in \RR$ and $C(u) > 0$ such that $(\theta + \epsilon u)^\top x + \varphi(x) \le D(u) - C(u) \|x\|$ for all $x \in \Rd$. By compactness, there exist finitely many unit vectors $u_1$, $u_2$, \ldots, $u_m$ such that the corresponding closed balls $B \bigl( u_i, (2\epsilon)^{-1} C(u_i) \bigr)$ cover the whole unit sphere in $\Rd$. Consequently, for any $x \in \Rd \setminus \{0\}$ and its direction $u(x) := \|x\|^{-1} x$, there exists an index $j = j(x) \in \{1,\ldots,m\}$ such that $\|u(x) - u_j\| \le (2\epsilon)^{-1} C(u_j)$, whence
\begin{eqnarray*}
	\theta^\top x + \epsilon \|x\| + \varphi(x)
	& = & (\theta + \epsilon u(x))^\top x + \varphi(x) \\
	& \le & (\theta + \epsilon u_j)^\top x + \varphi(x) + \epsilon \|u_j - u(x)\| \|x\| \\
	& \le & D(u_j) + \bigl( \epsilon \|u_j - u(x)\| - C(u_j) \bigr) \|x\| \\
	& \le & \max_{i=1,\ldots,m} D(u_i) - 2^{-1} \min_{i=1,\ldots,m} C(u_i) \|x\| \\
	& \to & - \infty	\quad\text{as} \ \|x\| \to \infty .
\end{eqnarray*}\\[-7ex]
\strut	\hfill	$\Box$

\paragraph{Proof of Theorem~\ref{thm:aebe}.}
As mentioned already, the statements about $\theta \in \Theta(P)$ and $\Pi(\cdot)$ are a consequence of Theorem~\ref{thm:dabigone}~(ii) and Proposition~\ref{prop:godot}. Note also that for $\theta \in \Rd \setminus \Theta(P)$ and arbitrary $r > 0$,
\begin{eqnarray*}
	\liminf_{n\to \infty} \int_{\Rd} \exp(\theta^\top x) \, P_n(dx)
	& \ge & 	\lim_{n\to \infty} \int_{\Rd} \min \bigl( \exp(\theta^\top x), r \bigr)
		\, P_n(dx) \\
	& = & \int_{\Rd} \min \bigl( \exp(\theta^\top x), r \bigr)
		\, P(dx) ,
\end{eqnarray*}
and the right hand side tends to infinity as $r \uparrow \infty$.
\hfill	$\Box$

\end{document}